\documentclass[11pt]{amsart}
\usepackage{amsmath,amssymb}
\usepackage{amscd,eucal,amsthm}
\newtheorem{Theorem}{Theorem}
\newtheorem{Corollary}{Corollary}
\newtheorem{Lemma}{Lemma}

\newtheorem{Proposition}{Proposition}

\theoremstyle{remark}

\renewcommand{\to}[1][]{\xrightarrow{#1}}

\newcommand{\C}{{\mathbb{C}}}

\newcommand{\B}{{\mathcal{B}}}

\newcommand{\Z}{{\mathbb{Z}}}
\renewcommand{\gg}{{\mathfrak{g}}}

\newcommand{\nn}{{\mathfrak{n}}}

\newcommand{\z}{{\mathfrak{z}}}
\renewcommand{\a}{{\mathfrak{a}}}

\newcommand{\bb}{{\mathfrak{b}}}
\newcommand{\hh}{{\mathfrak{h}}}
\renewcommand{\l}{{\lambda}}

\newcommand{\cc}{\mathbb{C}}

\newcommand{\VV}{\mathbb{V}}
\newcommand{\A}{\mathcal{A}}
\newcommand{\ZZ}{\mathcal{Z}} \DeclareMathOperator{\Fun}{Fun}
 \DeclareMathOperator{\Op}{Op}
 
 \DeclareMathOperator{\gr}{gr}
\DeclareMathOperator{\ad}{ad} \DeclareMathOperator{\ev}{ev}
 
 \DeclareMathOperator{\id}{id}

\DeclareMathOperator{\End}{End} \DeclareMathOperator{\Ind}{Ind}
\DeclareMathOperator{\Res}{Res} 
\DeclareMathOperator{\Spec}{Spec}
 \DeclareMathOperator{\rk}{rk}
\DeclareMathOperator{\Tr}{Tr} 
\renewcommand{\phi}{\varphi}

\newcommand{\ka}{\kappa}

 \makeatletter
\def\@mult#1{\raise #1\rlap{$\cdot$}\lower #1\rlap{$\cdot$}\cdot}
\def\did{\mathrel{\@mult{3pt}}}
\makeatother
\def\openrow#1#2#3{\setbox0=\vbox{\hbox
    {\vrule height#2 width#3\kern#2\vrule height#2 width0pt}\hrule height#3}
    \hbox{\leaders\copy0\hskip#1\wd0\vrule width#3}}
\def\row#1#2#3{\vbox{\hrule height#3\openrow{#1}{#2}{#3}}}
\def\Yr#1{\row{#1}{1.5ex}{.1ex}}

\def\DY#1\endDY{\baselineskip=1ex\lineskip=0pt\lineskiplimit=0pt{\vcenter
    {\Yr#1}}}
\def\openclm#1#2#3{\setbox0=\vbox{\hrule height#3\hbox
    {\vrule width0pt\kern#2\vrule width#3 height#2}}\vtop
    {\leaders\copy0\vskip#1\ht0\hrule height#3}}
\def\clm#1#2#3{\hbox{\vrule width#3\openclm{#1}{#2}{#3}}}
\def\Yc#1{\clm{#1}{1.5ex}{.1ex}}

\def\CDY#1\endCDY{{\vcenter{\hbox{\Yc#1}}}}

\newcommand{\nc}{\newcommand}
\nc{\g}{{\mathfrak g}} \nc{\LG}{{}^L\neg G} \nc{\pone}{{\mathbb
P}^1} \nc{\wt}{\widetilde} \nc{\wh}{\widehat} \nc{\ghat}{\wh{\gg}}
\nc{\mc}{\mathcal} \nc{\su}{{\mathfrak s}{\mathfrak l}}
\nc{\ppart}{(\!(t)\!)} \nc{\on}{\operatorname} \nc{\sw}{{\mf
s}{\mf l}} \nc{\mf}{\mathfrak} \nc{\ol}{\overline}
\nc{\Gr}{\on{Gr}} \nc{\bi}{\bibitem}

\nc{\reg}{\on{reg}}
\nc{\mb}{\mathbf}
\nc{\zpart}{(\!(z)\!)}

\begin{document}

\author[Boris Feigin]{Boris Feigin$^1$}

\thanks{$^1$ Supported by the grants RFBR 08-01-00720, RFBR
05-01-02934 and NSh-6358.2006.2.}

\address{Landau Institute for Theoretical Physics, Kosygina St 2,
Moscow 117940, Russia}

\author[Edward Frenkel]{Edward Frenkel$^2$}

\thanks{$^2$ Supported by DARPA and AFOSR through the grant
FA9550-07-1-0543.}

\address{Department of Mathematics, University of California,
  Berkeley, CA 94720, USA}

\author[Leonid Rybnikov]{Leonid Rybnikov$^3$}

\thanks{$^3$ Supported by
RFBR 05-01-02805-CNRSL-a, RFBR 07-01-92214-CNRSL-a, NSF grant
DMS-0635607 and Deligne 2004 Balzan prize in mathematics. L.R.
gratefully acknowledges the Institute for Advanced Study for
providing warm hospitality and excellent working conditions.}

\address{Institute for Theoretical and Experimental Physics,
B. Cheremushkinskaya 25, Moscow 117259, Russia}

\title{On the endomorphisms of Weyl modules over affine Kac-Moody
algebras at the critical level}

\date{February 2008}

\begin{abstract}
We present an independent short proof of the main result of
\cite{FG07} that the algebra of endomorphisms of a Weyl module of
critical level is isomorphic to the algebra of functions on the space
of monodromy-free opers on the disc with regular singularity and
residue determined by the highest weight of the Weyl module.  We
derive this from the results of \cite{FFR} about the shift of argument
subalgebras.
\end{abstract}

\maketitle

\section{Formulation of the main result}

\subsection{Weyl modules at the critical level.} Let $\gg$ be a simple
Lie algebra, and $\wh{\gg}$ be the corresponding affine Kac--Moody
algebra. The Lie algebra $\wh{\gg}$ is a central extension of the
formal loop algebra $\gg\ppart$ by one-dimensional center with
generator ${\mb 1}$. The commutation relations are as follows:
\begin{equation}
[g_1\otimes x(t),g_2\otimes y(t)]=[g_1,g_2]\otimes
x(t)y(t)+\kappa_c(g_1,g_2)\Res_{t=0}x(t)dy(t)\cdot {\mb 1},
\end{equation}
where $\kappa_c$ is the invariant inner product on $\gg$ defined
by the formula
\begin{equation}\kappa_c(g_1,g_2)=-\frac{1}{2}\Tr_{\gg}\ad(g_1)\ad(g_2).
\end{equation}
Set $\wh{\gg}_+=\gg[[t]]\subset\wh{\gg}$ and
$\wh{\gg}_-=t^{-1}\gg[t^{-1}]\subset\wh{\gg}$.

Define the completion $\wt{U}(\wh{\gg})$ of $U(\wh{\gg})$ as the
inverse limit of $U(\wh{\gg})/U(\wh{\gg})(t^n\gg[[t]])$, $n>0$.  The
action of $\wt{U}(\wh{\gg})$ is well-defined on the category of
discrete $\wh{\gg}$-modules, i.e., those in which every vector is
annihilated by $t^n\gg[[t]]$ for some $n>0$. We set
$$
\wt{U}_{\ka_c}(\wh{\gg})=\wt{U}(\wh{\gg})/({\mb 1}-1).
$$
This algebra acts on discrete $\wh{\gg}$-modules of {\em
critical level} (i.e., $\wh{\gg}$-modules on which the element $K$
acts as unity).

For a dominant integral weight $\l$ of $\gg$, let
$$\pi_\l:\gg\to\End_\cc V_\l$$ be the finite-dimensional irreducible
representation of $\gg$ with the highest weight $\l$.  One can
naturally extend this representation to $\wh{\gg}_+=\gg[[t]]$ by using
the composition with the natural map $\gg[[t]] \to \gg$ corresponding
to evaluation at $t=0$. The \emph{Weyl module} at the critical level
with the highest weight $\lambda$ is by definition the induced module
$$
\VV_\lambda:=\Ind_{\wh{\gg}_+ \oplus \C{\mb 1}}^{\wh{\gg}}V_\l,
$$
where ${\mb 1}$ acts on $V_\l$ as the identity.

\subsection{Action of the center and monodromy-free opers.}

Consider the Langlands dual Lie algebra ${}^L \gg$ whose Cartan matrix
is the transpose of the Cartan matrix of $\gg$. Denote by $^L G$ the
group of inner automorphisms of $^L \gg$. In \cite{FF,F:wak} the
center $Z(\ghat)$ of the completed enveloping algebra
$\wt{U}_{\ka_c}(\ghat)$ at the critical level was identified with the
algebra of polynomial functions on the space $\Op_{^L G}(D^\times)$ of
$^L G$-opers on the disc $D = \Spec \cc\ppart$.

Let us recall the notion of opers which was introduced in
\cite{BD}. Fix a Cartan decomposition
$$
^L \gg={}^L \nn\oplus {}^L \hh\oplus {}^L \nn_-.
$$
The Cartan subalgebra ${}^L \hh$ is canonically identified with
$\hh^*$. We denote by $\Pi^\vee$ the set of simple roots of $^L
\gg$ (which is the set of simple coroots of $\gg$). Set
$$
p_{-1} = \sum_{\alpha^\vee\in\Pi^\vee} e_{-\alpha^\vee}\in {}^L
\gg,
$$
where the $e_{-\alpha^\vee}$ are non-zero generators of the
$-\alpha^\vee$-root subspaces in $^L \nn_-$.

The space $\Op_{^L G}(D^\times)$ of $^L G$-opers is the quotient
of the space of connections on the trivial $^L G$-bundle on
$D^\times$ of the form
$$
d + (p_{-1} + {\mathbf v}(z))dz, \qquad {\mathbf v}(z) \in
{}^L\bb\zpart
$$
by the action of the group $^L N\zpart$.

Consider the action of the center $Z(\wh{\gg})$ on $\VV_\l$. Since
$Z(\wh{\gg})=\Fun(\Op_{^L G}(D^\times))$, the support of $\VV_\l$ as a
$Z(\wh{\gg})$-module is a closed subset in the space of opers $\Op_{^L
G}(D^\times)$.

In \cite{FG}, Sect. 2.9, a closed subspace $\Op_{^L
G}^{\l,\reg} \subset \Op_{^L G}(D^\times)$ of monodromy-free
opers with regular singularity and residue determined by $\l$ was
defined (this definition is reviewed in \cite{book}, Sect. 9.2.3, and
in \cite{FFTL}, Sect. 4.4).

The following assertion was proved in \cite{FG05}, Lemma 1.7.

\begin{Proposition}    \label{contained}
The support of $\VV_\l$ is contained in the subspace $\Op_{^L
G}^{\l,\reg} \subset \Op_{^L G}(D^\times)$.
\end{Proposition}

Furthermore, in \cite{FG07} the following theorem was proved which
completely describes the algebra of endomorphisms of the Weyl module
$\VV_\l$:

\begin{Theorem}    \label{main}
There is a commutative diagram
$$
\CD
Z(\ghat) @>{\sim}>> \Fun(\Op_{^L G}(D^\times)) \\
@VVV   @VVV  \\
\End_{\wh{\gg}}(\VV_\l) @>{\sim}>> \Fun(\Op_{^L
G}^{\l,\reg})
\endCD
$$
\end{Theorem}

The proof of this theorem given in \cite{FG07} used non-trivial
results about the semi-infinite cohomology of $\VV_\l$. The goal of
this paper is to give an alternative proof of this theorem, in which
we will not use semi-infinite cohomology, but will rely instead the
results about the shift of argument subalgebra from \cite{FFR} and on
Proposition \ref{contained}.

\subsection{Idea of the proof.} The loop rotation operator
$-t\partial_t$ acting on $\ghat$ defines a $\Z$-grading on
$Z(\wh{\gg})$ and on $\Fun(\Op_{^L G}(D^\times))$, and the isomorphism
of \cite{FF} preserves these gradings.

According to Proposition \ref{contained}, the action of $Z(\ghat)$ on
$\VV_\l$ factors through the algebra $\Fun(\Op_{^L
G}^{\l,\reg})$. Hence it is sufficient to prove that the map
$Z(\ghat)\to\End_{\wh{\gg}}(\VV_\l)$ is surjective, and that the
character of $\Fun(\Op_{^L G}^{\l,\reg})$ (understood as the formal
power series in a variable $q$ whose $q^n$ coefficient is the
dimension of the degree $n$ subspace) is not greater, term by term,
than that of $\End_{\wh{\gg}}(\VV_\l)$.

The character of $\Fun(\Op_{^L
G}^{\l,\reg})$ was computed in \cite{FG07}, Sect. 5.1, and is given by
the formula
\begin{equation}    \label{char opers}
\frac{\prod\limits_{\alpha>0}(1-q^{\langle\alpha^\vee,\l+\rho\rangle})}
{\prod\limits_{k=1}^\infty(1-q^k)^{\rk\gg}}.
\end{equation}
Here is a brief derivation of this formula. Recall that the algebra
$\Fun(\Op_{^L G}^{\l,\reg})$ is the quotient of the algebra
$\Fun(\Op_{^L G}(D^\times)^{RS}_\l)$, where $\Op_{^L G}^{\l,\on{RS}}$
is the space of opers with regular singularity and residue $\l$. The
latter is a free polynomial algebra with homogeneous generators
$P_i^{(k)}$, where $i=1,\dots,\rk\gg$, $k \geq 0$, whose degree is
equal to $k$. The subset $\Op_{^L G}^{\l,\reg}\subset\Op_{^L
G}^{\l,\on{RS}}$ is defined by a regular sequence of homogeneous
relations, enumerated by positive roots $\alpha$ of $\gg$, of the
degree $\langle\alpha^\vee,\l+\rho\rangle$ (see the above references
\cite{FG,book,FFTL} for details). This gives us formula (\ref{char
opers}).

In order to estimate the character of $\End_{\wh{\gg}}(\VV_\l)$, we
pass from $\End_{\wh{\gg}}(\VV_\l)$ to its associated graded algebra
with respect to the PBW filtration.  The crucial points in our proof
are the result of \cite{FFR} that $V_\l$ is a cyclic module over the
nilpotent shift of argument subalgebra $\A_f$ (which turns out to be
closely related to the associated graded algebra
$\gr\End_{\wh{\gg}}(\VV_\l)$) and the computation of invariants from
\cite{book}. Using these results, we obtain the desired lower bound
for the character of $\gr\End_{\wh{\gg}}(\VV_\l)$.

\subsection{Acknowledgements} We thank D. Gaitsgory for useful comments on the draft of this paper.

\section{Proof of the Theorem}

\subsection{Shift of argument subalgebras.}

To any $\mu\in\gg^*$ one can assign a commutative subalgebra
$\A_\mu\subset U(\gg)$ called the quantum shift of argument
subalgebra. This algebra comes from the center $Z(\ghat)$ in the
following way. Let $\pi:\wh{\gg}_+\to\gg$ be the homomorphism of
evaluation at $t=0$. Consider the following quantum Hamiltonian
reduction algebra
$$(U(\gg)\otimes\wt{U}_{\ka_c}(\wh{\gg})/U(\gg)\otimes
\wt{U}_{\ka_c}(\wh{\gg})
(\wh{\gg}_+-\pi(\wh{\gg}_+)))^{\wh{\gg}_+}.$$ Here
$U(\gg)\otimes\wt{U}_{\ka_c}(\wh{\gg})
(\wh{\gg}_+-\pi(\wh{\gg}_+))$ is the left ideal generated by
$x-\pi(x)$ for all $x\in\wh{\gg}_+$. The center $Z(\ghat)$ of
$\wt{U}_{\ka_c}(\ghat)$ naturally maps to this quotient. Each
element of the above quotient has a unique representative in
$U(\gg)\otimes U(\wh{\gg}_-)\subset
U(\gg)\otimes\wt{U}_{\ka_c}(\wh{\gg})$. Thus we obtain a
homomorphism $Z(\ghat)\to U(\gg)\otimes U(\wh{\gg}_-)$.

The element $\mu\in\gg^*$ defines a character
$\mu_-:\wh{\gg}_-\to\cc$ by $\mu_-(xt^{-1})=\langle\mu,x\rangle$
and $\mu_-(xt^{-k})=0$ for $k>1$. This gives us a homomorphism
$$
\id\otimes\mu_-:U(\gg)\otimes U(\wh{\gg}_-)\to U(\gg).
$$
Thus, we have a homomorphism $Z(\ghat)\to U(\gg)$ depending on
$\mu\in\gg^*$. The subalgebra $\A_\mu\subset U(\gg)$ is, by
definition, the image of this homomorphism (see \cite{R,FFTL} for more
details).

\medskip

Let ${\mf Z}$ be the image of the center $Z(\ghat)$ in $U(\gg)\otimes
U(\ghat_-)$. It follows from \cite{FFTL}, Theorem 5.6(1) and Lemma
5.5, that ${\mf Z}$ is a free commutative algebra generated by the
homogeneous (with respect to the grading defined by the loop-rotation
operator $-t\partial_t$) elements $P_i^{(k)}$ of degree $k$,
$i=1,\dots,\ell=\rk\gg,\ k=0,1,\dots$, such that $P_i^{(0)}\in
U(\gg)\otimes1$ are the generators of the center of $U(\gg)$. The
degree of $P_i^{(k)}$ with respect to the PBW filtration is $d_i+1$,
where $d_1,\dots,d_\ell$ are the exponents of $\gg$. From the
description of these elements given in the above reference it is easy
to see that the associated graded of ${\mf Z}$ with respect to the PBW
filtration on the second factor of $U(\gg)\otimes U(\ghat_-)$ is
freely generated by elements $\overline{P_i^{(k)}}, i=1,\dots,\rk\gg,\
k\geq 0$, such that $\overline{P_i^{(k)}}\in U(\gg)\otimes S(\ghat_-)$
for $i=1,\dots,\rk\gg,\ k=0,1,\dots, d_i$, where the $d_i$ are the
exponents of $\gg$, and $\overline{P_i^{(k)}}\in 1\otimes S(\ghat_-)$
for $i=1,\dots,\rk\gg,\ k>d_i$.

\begin{Lemma}    \label{images} (see also
\cite{FFTL}, Lemma 3.13)
For regular $\mu$ the images of ${\mf Z} \subset U(\g) \otimes
U(\ghat_-)$ and $\gr {\mf Z} \subset U(\g) \otimes S(\ghat_-)$ under
$\id \otimes \mu_-$ and $\id \otimes \gr\mu_-$, respectively, coincide
and are equal to the same commutative subalgebra ${\mc A}_\mu \subset
U(\gg)$.
\end{Lemma}

\begin{proof}
According to \cite{R}, Theorem 1, and \cite{FFTL}, Theorem 5.8, for
any regular $\mu\in\gg^*$ the subalgebra $\A_\mu\subset U(\gg)$ is
freely generated by the images of $P_i^{(k)}$ with $i=1,\dots,\rk\gg,\
k=0,1,\dots, d_i$, where the $d_i$ are the exponents of $\gg$. Since
the elements $P_i^{(k)}$ are homogeneous with respect to the loop
rotation operator, the image of $P_i^{(k)}$ in $U(\g) \otimes S(\gg
t^{-1})=U(\g) \otimes U(t^{-1}\gg[t^{-1}]/(1\otimes
t^{-2}\gg[t^{-1}])$ is homogeneous with respect to the grading on the
second factor, and hence coincides with the image of
$\overline{P_i^{(k)}}$ in $U(\g) \otimes S(\gg t^{-1})=U(\g) \otimes
S(t^{-1}\gg[t^{-1}]/(1\otimes t^{-2}\gg[t^{-1}])$. Since the
homomorphism $\id \otimes \mu_-$ factors through $U(\g) \otimes S(\gg
t^{-1})$, we have
$$
\id\otimes\mu_-(P_i^{(k)})=\id\otimes\gr\mu_-(\overline{P_i^{(k)}}).
$$
Therefore the images of $\overline{P_i^{(k)}}\in U(\gg)\otimes
S(\ghat_-)$ with $k=0,1,\dots,d_i$ under the map
$\id\otimes\gr\mu_-$ generate the same commutative subalgebra $\A_\mu
\subset U(\gg)$.
\end{proof}

We remark that a certain limit of ${\mc A}_\mu$ in the case when
$\g=\sw_n$ may be identified with the Gelfand--Zetlin algebra (see
\cite{R}). Hence the algebra ${\mc A}_\mu$ may be thought of as a
generalization of the Gelfand--Zetlin algebra to an arbitrary simple
Lie algebra.

\bigskip

An important special case is when $\mu \in \gg^* \simeq \gg$ is a
regular nilpotent element. Since all of these elements belong to a
single coadjoint orbit, it is sufficient to consider one particular
representative. Let
\begin{equation}    \label{f}
f=\sum\limits_{\alpha\in\Pi} e_{-\alpha}\in\gg \simeq \gg^*
\end{equation}
(the last isomorphism is obtained from any non-degenerate inner
product on $\gg$ which we fix once and for all) be the principal
nilpotent element. Let $\{ e,\ h,\ f \}$ be a principal $\sw_2$-triple
in $\gg$ containing $f$. The operator $\ad h$ defines a gradation on
$U(\gg)$ which is called the \emph{principal gradation}. The algebra
$\A_f$ is generated by homogeneous elements with respect to the
principal gradation on $U(\gg)$. Moreover, the homomorphism
$Z(\ghat)\to \A_f$ is a homomorphism of \emph{graded} algebras. The
algebra $\A_f$ acts on $V_\l$ by creation operators.

The following result was proved in \cite{FFR}:

\begin{Theorem}    \label{cyclic}
The module $V_\l$ is cyclic as an $\A_\mu$-module for any regular $\mu
\in \g^*$. Moreover, if $\mu = f$, then the highest weight vector of
$V_\l$ is a cyclic vector.
\end{Theorem}

\subsection{The associated graded of $\End_{\wh{\gg}}(\VV_\l)$.}

Each endomorphism of the Weyl module $\VV_\l$ is uniquely determined
by the image of the generating subspace $V_\l \subset \VV_\l$. Hence
the algebra $\End_{\wh{\gg}}(\VV_\l)$ may be naturally identified with
$$\left(\End_\cc(V_\l)\otimes U_{\ka_c}(\wh{\gg})/\End_\cc(V_\l)\otimes
U_{\ka_c}(\wh{\gg})(\wh{\gg}_+-\pi_\l(\wh{\gg}_+))\right)^{\wh{\gg}_+}.$$
Here $\End_\cc(V_\l)\otimes
U_{\ka_c}(\wh{\gg})(\wh{\gg}_+-\pi_\l(\wh{\gg}_+))$ is the left ideal
generated by $x-\pi_\l(x)$ for all $x\in\wh{\gg}_+$. Since each
element of the above quotient algebra has a unique representative
in $\End_\cc(V_\l)\otimes U(\wh{\gg}_-)$, the algebra
$\End_{\wh{\gg}}(\VV_\l)$ may be regarded as a subalgebra
$\B\subset\End_\cc(V_\l)\otimes U(\wh{\gg}_-)$ (in the same way as
in the previous subsection).

Let $\ZZ\subset\B=\End_{\wh{\gg}}(\VV_\l)$ be the image of the center
of the completed enveloping algebra at the critical level $Z(\ghat)
\subset \wt{U}_{\ka_c}(\wh{\gg})$. Since each element of $\B$ commutes
with $\wt{U}_{\ka_c}(\wh{\gg})$, we find that $\ZZ$ belongs to the
center of $\B$.

Consider a filtration on the algebra $\End_\cc(V_\l)\otimes
U(\wh{\gg}_-)$ determined by the trivial filtration on the first
factor and the PBW filtration on the second one. This filtration
determines (by restriction to a subalgebra) filtrations on $\B$ and
$\ZZ$. Due to the PBW theorem, the associated graded of
$\End_\cc(V_\l)\otimes U(\wh{\gg}_-)$ is
$$\End_\cc(V_\l)\otimes S(\wh{\gg}_-) \simeq \End_\cc(V_\l)\otimes
\Fun(\gg[[t]]).$$ The isomorphism $S(\wh{\gg}_-) \simeq
\Fun(\gg[[t]])$ depends on the choice of an invariant non-degenerate
scalar product $(\cdot, \cdot)$ on $\gg$. Namely, a linear element
$xt^{-k}\in S(\wh{\gg}_-)$ maps to a linear function defined on
$y(t)\in\gg[[t]]$ as $\Res_{t=0}(x,y(t))t^{-k}$. In particular, $S(\gg
t^{-1})$ does to $\Fun(\gg)$ under this isomorphism.

The associated graded of $\B$,
$$\overline{\B}:=\gr\B=\gr\left(\End_\cc(V_\l)\otimes
U_{\ka_c}(\wh{\gg})/\End_\cc(V_\l)\otimes
U_{\ka_c}(\wh{\gg})(\wh{\gg}_+-\pi_\l(\wh{\gg}_+))\right)^{\wh{\gg}_+},$$
is naturally embedded into
\begin{multline}    \label{image}
\left(\gr(\End_\cc(V_\l)\otimes
U_{\ka_c}(\wh{\gg})/\End_\cc(V_\l)\otimes
U_{\ka_c}(\wh{\gg})(\wh{\gg}_+-\pi_\l(\wh{\gg}_+)))\right)^{\wh{\gg}_+}
\\ \simeq \left(\End_\cc(V_\l)\otimes
\Fun(\gg[[t]])\right)^{\gg[[t]]}.
\end{multline}

However, it was shown in \cite{book} that this embedding is not an
isomorphism unless $\l=0$ or minuscule (the reason for this is that it
is only for these $\l$ that the module $V_\l$ is cyclic for the
centralizer $\a_f$ of the principal nilpotent element $f$).
Nevertheless, we will now use our results on the shift of argument
subalgebra from \cite{FFR} to give an estimate of the image of
$\overline{\B}$ in \eqref{image}, which will turn out to be sufficient
for our purposes.

The subalgebra
$$
\overline{\ZZ}:=\gr\ZZ\subset\left(\End_\cc(V_\l)\otimes
\Fun(\gg[[t]])\right)^{\gg[[t]]}
$$
is generated by the elements $\pi_\l\otimes1(\overline{P_i^{(k)}})$
with $i=1,\dots,\rk\gg,\ k \geq 0$ (this follows from the definition
of the elements $\overline{P_i^{(k)}}$).

Consider the subalgebra
$$
\overline{\ZZ'}:=1\otimes
\Fun(\gg[[t]])^{\gg[[t]]}\subset\left(\End_\cc(V_\l)\otimes
\Fun(\gg[[t]])\right)^{\gg[[t]]}.
$$

\begin{Lemma} \label{Z' in Z}
The algebra $\overline{\ZZ'}$ is a free polynomial algebra generated
by $\pi_\l\otimes1(\overline{P_i^{(k)}})$ with $i=1,\dots,\rk\gg,\ k
\geq d_i+1$, and hence $\overline{\ZZ'}\subset\overline{\ZZ}$.
\end{Lemma}

\begin{proof}
According to a result of \cite{BD} (see \cite{book}, Theorem 3.4.2,
for a proof), the algebra $\overline{\ZZ'}$ is a free polynomial
algebra generated by some homogeneous elements $\overline{S_i^{(k)}}$
with $i=1,\dots,\rk\gg,\ k \geq d_i+1$ of degrees $k$ with respect to
the loop-rotation grading and $d_i+1$ with respect to the grading by
the degree of polynomials. Hence the elements $\overline{P_i^{(k)}}$
with $i=1,\dots,\rk\gg,\ k \geq d_i+1$ generate a free polynomial
subalgebra of the same size in $\left(U(\gg)\otimes
\Fun(\gg[[t]])\right)^{\gg[[t]]}$. Thus it remains to show that
$\overline{P_i^{(k)}}\in1\otimes \Fun(\gg[[t]])^{\gg[[t]]}$. Note that
the elements $P_i^{(k)}\in\ZZ$ with $i=1,\dots,\rk\gg,\ k \geq d_i+1$
are homogeneous with respect to the loop-rotation grading, and hence
their leading terms with respect to the PBW filtration belong to
$1\otimes U(\ghat_-)$. Therefore $\overline{P_i^{(k)}}\in1\otimes
\Fun(\gg[[t]])^{\gg[[t]]}$. Hence the assertion.
\end{proof}

Let $\mathcal{I}$ be the left ideal in $\left(\End_\cc(V_\l)\otimes
\Fun(\gg[[t]])\right)^{\gg[[t]]}$ generated by all
$\overline{P_i^{(k)}}$ ($i=1,\dots,\rk\gg,\ k \geq d_i+1$).  According
to \cite{book}, Sects. 9.6.4--9.6.5, the algebra
$\left(\End_\cc(V_\l)\otimes \Fun(\gg[[t]])\right)^{\gg[[t]]}$ is a
free $\overline{\ZZ'}$-module. Any space of generators of this module
is therefore isomorphic to the space of $\mathcal{I}$-coinvariants
$$\left(\End_\cc(V_\l)\otimes
\Fun(\gg[[t]])\right)^{\gg[[t]]}/\mathcal{I}
\simeq\left(\End_\cc(V_\l)\otimes
\Fun(\gg)\right)^{\gg}/1\otimes\Fun(\gg)^{\gg}_+,
$$
where
$\Fun(\gg)^{\gg}_+$ is the maximal graded
ideal in $\Fun(\gg)^{\gg}$.

According to \cite{K}, the latter quotient has the following
description. Let $f\in\gg=\gg\cdot1\subset\gg[[t]]$ be the principal
nilpotent element \eqref{f}. The evaluation homomorphism at $f$,
$$\id\otimes\ev_f:\End_\cc(V_\l)\otimes
\Fun(\gg[[t]])\to\End_\cc(V_\l)$$ annihilates $\mathcal{I}$ and
gives rise to an isomorphism
$$\left(\End_\cc(V_\l)\otimes
\Fun(\gg[[t]])\right)^{\gg[[t]]}/\mathcal{I}\simeq
\End_{\a_f}(V_\l),$$ where $\a_f\subset\gg$ is the centralizer of
$f$. This is an isomorphism of \emph{graded} algebras with respect to
the loop-rotation grading (defined by the operator $-t\partial_t$) on
the left-hand side and the principal grading on the right-hand side.

Let now $\A_f\subset U(\gg)$ be the quantum shift of argument
subalgebra corresponding to the principal nilpotent element
$f\in\gg=\gg^*$.

\begin{Lemma}\label{diagram} There is a commutative diagram of graded
  algebras (with respect to the loop-rotation grading on the left-hand
  side and the principal grading on the right-hand side).
$$ \CD \overline{\ZZ}/(\mathcal{I}\cap\overline{\ZZ}) @>{\sim}>>
\pi_\l(\A_f) \\ @V{\sim}VV @V{\sim}VV \\
\overline{\B}/(\mathcal{I}\cap\overline{\B}) @>{\sim}>>
\End_{\A_f}(V_\l) \\ @VVV @VVV \\ \left(\End_\cc(V_\l)\otimes
\Fun(\gg[[t]])\right)^{\gg[[t]]}/\mathcal{I} @>{\sim}>>
\End_{\a_f}(V_\l)
\endCD
$$
\end{Lemma}

\begin{proof}
By definition of the quantum shift of argument subalgebra,
$\id\otimes f_-(\ZZ)$ is $\pi_\l(\A_f)\subset\End_\cc(V_\l)$.
Since $\gr f_-=\ev_f$, Lemma \ref{images} implies that
$\id\otimes\ev_f(\overline\ZZ)=\pi_\l(\A_f)$. Since $\ZZ$ belongs
to the center of $\B$, we have
$\id\otimes\ev_f(\overline{\B})\subset \End_{\A_f}(V_\l)$.

On the other hand, according to Theorem \ref{cyclic}, the algebra $\A_f$
has a cyclic vector in $V_\l$. Hence
$$\End_{\A_f}(V_\l)=\pi_\l(\A_f) =
\id\otimes\ev_f(\overline\ZZ)\subset\id\otimes\ev_f(\overline{\B}).$$
Thus, we obtain the opposite inclusion, which proves the assertion of
the lemma.
\end{proof}

Now we are going to prove that $\overline{\B} = \overline{\ZZ}$. Note
that the algebra $\overline{\ZZ}$ is generated over $\overline\ZZ'$ by
the elements $\pi_\l\otimes1(\overline{P_i^{(k)}})$ with
$i=1,\dots,\rk\gg,\ k \leq d_i$. Due to the same "homogeneity"
argument as in Lemmas~\ref{images}~and~\ref{Z' in Z} we find that
$\pi_\l\otimes1(\overline{P_i^{(k)}})$ with $i=1,\dots,\rk\gg,\ k \leq
d_i$ belong to $\End_\cc(V_\l)\otimes S(\gg
t^{-1})=\End_\cc(V_\l)\otimes \Fun(\gg)$. Hence we can lift the graded
space $\overline{\ZZ}/(\mathcal{I}\cap\overline{\ZZ})=\overline{\B}/
(\mathcal{I}\cap\overline{\B})$ to a graded subspace
$N\subset\overline{\ZZ}\cap\Fun(\gg)$.

\begin{Lemma}\label{free modules} $\overline{\B}$ and $\overline{\ZZ}$
are both free $\overline\ZZ'$-modules generated by $N$.
\end{Lemma}

\begin{proof}
Let $M\subset\left(\End_\cc(V_\l)\otimes
\Fun(\gg)\right)^{\gg}$ be a space of generators of
$$\left(\End_\cc(V_\l)\otimes\Fun(\gg[[t]])\right)^{\gg[[t]]}$$ as a
free $\overline{\ZZ'}$-module, containing $N$.
Clearly, $M$ also
freely generates the $\Fun(\gg)^{\gg}$-module $
\left(\End_\cc(V_\l)\otimes \Fun(\gg)\right)^{\gg}.$

For each regular $\mu_0\in\gg$, the evaluation at $\mu_0$ gives us an
isomorphism of vector spaces
$$\id\otimes\ev_{\mu_0}:M\overset
\sim\longrightarrow\End_{\z_{\gg}(\mu_0)}V_\l,$$ where
$\z_{\gg}(\mu_0)$ is the centralizer of $\mu_0$ in $\gg$ (this is a
classical result due to Kostant \cite{K}). Since
$N\subset\overline{\ZZ}$, we have
$\id\otimes\ev_{\mu_0}(N)\subset\pi_\l(\A_{\mu_0})$. Since
$\id\otimes\ev_{\mu_0}$ is an injection on $M$, we have $\dim
(\id\otimes\ev_{\mu_0}(N))=\dim N = \dim
V_\l=\dim\pi_\l(\A_{\mu_0})$. Hence, for each regular $\mu_0\in\gg$,
the evaluation at $\mu_0$ gives us an isomorphism of vector spaces
$$
\id\otimes\ev_{\mu_0}:N \overset
\sim\longrightarrow\pi_\l(\A_{\mu_0}).
$$
Therefore, for any $\mu=\mu_0+t\mu_1+\dots\in\gg[[t]]$, we have
$\id\otimes\ev_\mu(\overline\ZZ' \cdot N)=\pi_\l(\A_{\mu_0})$.

The submodules $\overline{\B}\subset \overline\ZZ' \cdot M$ and
$\overline{\ZZ}\subset \overline\ZZ' \cdot M$ clearly contain
$\overline\ZZ' \cdot N$. Since $\overline\ZZ' \cdot
N\subset\overline\ZZ$, the subalgebras
$\id\otimes\ev_\mu(\overline\ZZ)\subset\End_\cc V_\l$ and
$\ev_\mu(\overline{\B})\subset\End_\cc V_\l$ belong to the commutant
of $\id\otimes\ev_\mu(N)=\pi_\l(\A_{\mu_0})$ in $\End_\cc V_\l$ for
any $\mu=\mu_0+t\mu_1+\dots\in\gg[[t]]$ with regular $\mu_0$. By
Theorem \ref{cyclic}, $V_\l$ is a cyclic $\A_{\mu_0}$-module and hence
$\pi_\l(\A_{\mu_0}) =\End_{\A_{\mu_0}}(V_\l)$. Thus we have
$$\id\otimes\ev_\mu(\overline\ZZ' \cdot N)=\pi_\l(\A_{\mu_0})
=\End_{\A_{\mu_0}}(V_\l)=
\id\otimes\ev_\mu(\overline{\B})=\id\otimes\ev_\mu(\overline\ZZ).$$
Since $\overline\ZZ'$ is the algebra of $\gg[[t]]$-invariant functions
on $\gg[[t]]$, each point of $\Spec \overline\ZZ'$ has a
representative in $\gg[[t]]$ of the form
$\mu=\mu_0+t\mu_1+\dots\in\gg[[t]]$ with regular $\mu_0$.  Therefore
we find that the images of $\overline{\B}$ and $\overline{\ZZ}$ in the
quotient of $\overline\ZZ' \cdot M$ by each maximal ideal of
$\overline\ZZ'$ have the same dimension equal to $\dim N=\dim V_\l$.

Let $C$ be a complementary subspace to $N$ in the space of generators
$M$, and let $c_1,...,c_n$ be a basis of $C$. Then $\overline\ZZ'
\cdot M=\overline\ZZ' \cdot N\oplus \overline\ZZ' \cdot C$. Suppose
that $\overline{\B}\ne \overline\ZZ' \cdot N$. Then the module
$\overline{\B}$ contains an element from $\overline\ZZ' \cdot C$ of
the form $z_1c_1+...+z_nc_n$ with $z_i\in \overline\ZZ'$. Choose a
maximal ideal $J\subset\overline\ZZ'$ such that $z_1\not\in J$ (i.e.
$\mu=\mu_0+t\mu_1+\dots\in\gg[[t]]$ with regular $\mu_0$ such that
$z_1(\mu)\ne0$). Then $\overline{\B}/(\overline{\B}\cap J \cdot
M)\supset N+c_1$ and therefore $\dim\overline{\B}/(\overline{\B}\cap J
\cdot M)>\dim N$. On the other hand, we have found above that
$\dim\overline{\B}/(\overline{\B}\cap J \cdot M)=\dim N$ for all
maximal ideals $J$. Hence we obtain a contradiction. In the same way
we prove that $\overline{\ZZ}=\overline{\ZZ} \cdot N$.
\end{proof}

\begin{Corollary}
$\B = \ZZ$.
\end{Corollary}

Thus, the homomorphism $Z(\ghat)\to\End_{\wh{\gg}}(\VV_\l)$ is
surjective. Hence it remains to show that the character of $\ZZ$ (or,
equivalently, $\overline{\ZZ}$) is not smaller than that of
$\Fun(\Op_{^L G}^{\l,\reg})$. This will be done in the next section.

\subsection{Comparison of characters}

According to \cite{FFR}, the character of $\pi_\l(\A_f)$ is
the same as that of $V_\l$ with respect to the principal grading,
which is known to be
$$
\frac{\prod\limits_{\alpha>0}(1-q^{\langle\alpha^\vee,\l+\rho\rangle})}
{\prod\limits_{\alpha>0}(1-q^{\langle\alpha^\vee,\rho\rangle})}.
$$
The denominator may be rewritten as
$$\prod\limits_{\alpha>0}(1-q^{\langle\alpha^\vee,\rho\rangle})=
\prod\limits_{i=1}^{\rk\gg}\prod\limits_{k=1}^{d_i}(1-q^k).
$$
On the other hand, the character of $\overline{\ZZ'}$ (with respect to the loop-rotation grading) is
$\prod\limits_{i=1}^{\rk\gg}\prod\limits_{k=d_i+1}^{\infty}(1-q^k)$.
According to Lemma~\ref{diagram} and Lemma~\ref{free modules}, $\overline{\ZZ}$ is
a free $\overline\ZZ'$-module with the space of generators
$\pi_\l(\A_f)$. Hence the character of
$\overline{\ZZ}$ is the product of those of
$\pi_\l(\A_f)$ and $\overline\ZZ'$. Therefore it is equal to
$$
\frac{\prod\limits_{\alpha>0}(1-q^{\langle\alpha^\vee,\l+\rho\rangle})}
{\prod\limits_{k=1}^\infty(1-q^k)^{\rk\gg}},
$$
which coincides with the character of $\Fun(\Op_{^L
G}^{\l,\reg})$ given by formula \eqref{char opers}.

Since $\ZZ$ factors through $\Fun(\Op_{^L G}^{\l,\reg})$ by
Proposition \ref{contained} and we have shown that the homomorphism
$Z(\ghat)\to\End_{\wh{\gg}}(\VV_\l)$ is surjective, this completes the
proof of Theorem \ref{main}.


\begin{thebibliography}{FFTL}
\bibitem[BD]{BD} A. Beilinson and V. Drinfeld, {\em Quantization of
Hitchin's integrable system and Hecke eigensheaves}, Preprint,
available at www.ma.utexas.edu/$\sim$benzvi/BD.

\bibitem[FF]{FF} B. Feigin and E. Frenkel, {\em Affine Kac--Moody
algebras at the critical level and Gelfand--Dikii algebras}, Int.
Jour. Mod. Phys. {\bf A7}, Supplement 1A (1992) 197--215.

\bibitem[FFR]{FFR} B. Feigin, E. Frenkel and L. Rybnikov, \emph{Opers with
irregular singularity and spectra of the shift of argument
subalgebra.} Preprint math.QA/0712.1183.

\bibitem[FFT]{FFTL} B. Feigin, E. Frenkel and V. Toledano Laredo, {\em
Gaudin model with irregular singularities}, Preprint
math.QA/0612798.

\bibitem[FG07]{FG07} E. Frenkel and D. Gaitsgory, \emph{Weyl modules
and opers without monodromy}, Preprint math.QA/0706.3725.

\bibitem[FG05a]{FG} E. Frenkel and D. Gaitsgory {\em Local geometric
Langlands correspondence and affine Kac--Moody algebras},
Algebraic Geometry and Number Theory, Progress in Math. {\bf 253},
pp. 69--260, Birkh\"auser Boston, 2006 (math.RT/0508382).

\bibitem[FG05b]{FG05} E. Frenkel and D. Gaitsgory, \emph{Fusion and
convolution: applications to affine Kac-Moody algebras at the critical
level}, Pure and Applied Math. Quart. {\bf 2}, (2006) 1255--1312
(math.RT/0511284).

\bibitem[Fr05]{F:wak} E.~Frenkel, {\em Wakimoto modules, opers and the
center at the critical level}, Adv. Math. {\bf 195} (2005) 297--404
(math.QA/0210029).

\bibitem[Fr07]{book} E. Frenkel, \emph{Langlands correspondence for
loop groups.}  Cambridge Studies in Advanced Mathematics,
103. Cambridge University Press, Cambridge, 2007.

\bibitem[K]{K} B. Kostant, \emph{Lie group representations on polynomial
rings}, Amer. J. Math.  {\bf 85} (1963) 327--404.

\bibitem[R]{R} L. Rybnikov, {\em Argument shift method
and Gaudin model}, Func. Anal. Appl. {\bf 40} (2006), No. 3,
translated from Funktsional'nyi Analiz i Ego Prilozheniya {\bf 40}
(2006), No.~3, pp.~30--43 (math.RT/0606380).



\end{thebibliography}
\end{document}